\documentstyle[leqno,11pt]{article}
\begin{document}

\title{\bf On exponentially coprime integers}
\author{{\sc L\'aszl\'o T\'oth} (P\'ecs, Hungary)}
\date{Pure Math. Appl. (PU.M.A.), 15 (2004), 343-348}
\maketitle

\begin{abstract} The integers $n=\prod_{i=1}^r p_i^{a_i}$ and $m=\prod_{i=1}^r p_i^{b_i}$
having the same prime factors are called exponentially coprime if $(a_i,b_i)=1$ for
every $1\le i\le r$. We estimate the number of pairs of exponentially coprime
integers $n,m\le x$ having the prime factors $p_1,...,p_r$ and show that the
asymptotic density of pairs of exponentially coprime integers having $r$ fixed
prime divisors is $(\zeta(2))^{-r}$. \end{abstract}

\centerline{Mathematics Subject Classification (2000): 11A05, 11A25, 11N37}

\vskip5mm

{\bf 1. Introduction} \vskip3mm

Let $n>1$ be an integer of canonical form $n=\prod_{i=1}^r p_i^{a_i}$.
The integer $d$ is called an {\sl exponential divisor} of $n$ if
$d=\prod_{i=1}^r p_i^{c_i}$, where $c_i | a_i$ for every $1\le i \le r$, notation: $d|_e n$.
By convention $1|_e 1$. This notion was introduced by {\sc M. V. Subbarao} \cite{Su72}.
The smallest exponential divisor of $n>1$ is its squarefree kernel
$\kappa(n):=\prod_{i=1}^r p_i$.

Let $\tau^{(e)}(n)= \sum_{d|_e n} 1$ and $\sigma^{(e)}(n)=\sum_{d|_e n} d$ denote the
number and the sum of exponential divisors of $n$, respectively. Properties of these
functions were investigated by several authors, see \cite{FaSu89},  \cite{KaSu2003},
\cite{PeWu97}, \cite{SmWu97}, \cite{Su72}, \cite{Wu95}.

Two integers $n,m >1$ have common exponential divisors iff they have the same prime
factors and for $n=\prod_{i=1}^r p_i^{a_i}$, $m=\prod_{i=1}^r p_i^{b_i}$, $a_i,b_i\ge 1$
($1\le i\le r$), the {\sl greatest common exponential divisor} of $n$ and $m$ is
$$
(n,m)_e:=\prod_{i=1}^r p_i^{(a_i,b_i)}.
$$

Here $(1,1)_e=1$ by convention and $(1,m)_e$ does not exist for $m>1$.

The integers $n,m >1$ are called {\sl exponentially coprime}, if they have the same
prime factors and $(a_i,b_i)=1$ for every $1\le i\le r$, with the notation of above.
In this case $(n,m)_e=\prod_{i=1}^r p_i$. $1$ and $1$ are considered to be
exponentially coprime. $1$ and $m>1$ are not exponentially coprime. Exponentially
coprime integers were introduced by {\sc J. S\'andor} \cite{Sa96}.

Let $p_i$ ($1\le i\le r$) be fixed distinct primes and let $P^{(e)}(p_1,...,p_r;x)$
denote the number of pairs $\langle n,m \rangle$ of exponentially coprime integers
such that $\kappa(n)=\kappa(m)=\prod_{i=1}^r p_i$ and $n,m\le x$.

In this note we estimate $P^{(e)}(p_1,...,p_r;x)$ and show that the asymptotic
density of pairs of exponentially coprime integers having $r$ fixed prime divisors is
$(\zeta(2))^{-r}$.

As an open problem we formulate the following: What can be said on the asymptotic density
of pairs of exponentially coprime integers if their prime divisors are not fixed ?

\vskip3mm
{\bf 2. Results}
\vskip3mm

For a real $x\ge 1$ and an integer $n\ge 1$ consider the Legendre-type function
$L^{(e)}(x,n)$ defined as the number of integers $k\le x$ such that $k$ and $n$ are
exponentially coprime.

The following estimate holds:

\vskip2mm {\bf Theorem 1.} {\it Let $r\ge 1$ be a fixed integer and $p_1,...,p_r$ be
fixed distinct primes. Then uniformly for the real $x\ge 3$ and $n=\prod_{i=1}^r
p_i^{a_i}$ with $a_1,...,a_r\ge 1$,
$$
L^{(e)}(x,n)= \frac1{r!} \left( \prod_{i=1}^r \frac{\phi(a_i)}{a_i \log p_i} \right)
(\log x)^r + O\left( (\log x)^{r-1} \sum_{i=1}^r \theta(a_i) \right), \leqno(1)
$$
where $\phi(a)$ is Euler's function and $\theta(a)$ denotes the number of squarefree divisors of $a$.}
\vskip2mm

Let $N(p_1,...,p_r;x)$ denote the number of integers $n\le x$ having the kernel
$\kappa(n)=p_1\cdots p_r$. Taking $a_1=\cdots =a_r=1$ we obtain from Theorem 1 the
following known estimate, cf. for ex. \cite{Te95}, Ch. III.5 regarding integers free
of large prime factors.

\vskip2mm {\bf Corollary 1.} {\it
$$ N(p_1,...,p_r;x)=
\frac1{r!} \left( \prod_{i=1}^r \frac{1}{\log p_i} \right)
(\log x)^r + O((\log x)^{r-1}), \leqno(2)
$$ }

\vskip2mm {\bf Theorem 2.} {\it Let $r\ge 1$ be a fixed integer and $p_1,...,p_r$ be
fixed distinct primes. Then
$$
P^{(e)}(p_1,...,p_r;x)
= \frac1{(r!)^2(\zeta(2))^r} \left( \prod_{i=1}^r \frac1{(\log p_i)^2}\right)
(\log x)^{2r} + O\left( (\log x)^{2r-1} \log \log x \right). \leqno(3)
$$ }

\vskip2mm {\bf Corollary 2.} {\it The asymptotic density of pairs of exponentially
coprime integers having $r$ fixed prime divisors is $(\zeta(2))^{-r}$.} \vskip2mm

\vskip3mm
{\bf 3. Proofs}
\vskip3mm

The proofs of Theorems 1 and 2 are by induction on $r$, while Corollary 3 follows
from Theorem 2 and Corollary 1. First we prove the following lemma.

\vskip2mm
{\bf Lemma 1.} {\it Let $r\ge 1$ be a fixed integer and $t_1,...,t_r>0$ be fixed
real numbers. Then uniformly for the real $z\ge 1$ and the integers $a_1,...,a_r\ge 1$,
$$
\sum_{k_1t_1+\cdots +k_rt_r\le z \atop{(k_1,a_1)=\cdots =(k_r,a_r)=1 \atop{k_1,...,k_r\ge 1}}} 1 = \frac1{r!} \left( \prod_{i=1}^r \frac{\phi(a_i)}{a_i t_i} \right) z^r +
O\left(z^{r-1} \sum_{i=1}^r \theta(a_i) \right). \leqno(4)
$$ }
\vskip2mm

\vskip2mm {\bf Proof of Lemma 1.} We will use the well-known estimate: if $s\ge 0$,
then
$$
\phi_s(z,a):=\sum_{n\le z \atop{(n,a)=1}} n^s =\frac{z^{s+1}\phi(a)}{(s+1)a}+O(z^s\theta(a)), \leqno(5)
$$
uniformly for $z\ge 1$ and $a\ge 1$.

Induction on $r$. For $r=1$ (4) follows from (5) applied for $s=0$. Suppose formula
(4) is valid for $r-1$ and prove it for $r$.
$$
\sum_{k_1t_1+\cdots +k_rt_r\le z \atop{(k_1,a_1)=\cdots =(k_r,a_r)=1 \atop{k_1,...,k_r\ge 1}}} 1 =
\sum_{k_rt_r\le z-t_1-\cdots -t_{r-1} \atop{(k_r,a_r)=1 \atop{k_r\ge 1}}}
\sum_{k_1t_1+\cdots +k_{r-1}t_{r-1}\le z-k_rt_r \atop{(k_1,a_1)=\cdots =(k_{r-1},a_{r-1})=1 \atop{k_1,...,k_{r-1}\ge 1}}} 1
$$
$$
=\sum_{k_rt_r\le z-t_1-\cdots -t_{r-1} \atop{(k_r,a_r)=1 \atop{k_r\ge 1}}}
\left(\frac1{(r-1)!} \left( \prod_{i=1}^{r-1} \frac{\phi(a_i)}{a_i t_i}\right)
(z-k_rt_r)^{r-1} + O\left( z^{r-2} \sum_{i=1}^{r-1} \theta(a_i) \right) \right)
$$
$$
= \frac1{(r-1)!} \left( \prod_{i=1}^{r-1} \frac{\phi(a_i)}{a_i t_i}\right)
\sum_{k_rt_r\le z-t_1-\cdots -t_{r-1} \atop{(k_r,a_r)=1 \atop{k_r\ge 1}}} (z-k_rt_r)^{r-1}+
O\left( z^{r-1} \sum_{i=1}^{r-1} \theta(a_i) \right).
$$
Using the binomial formula and estimate (5) the sum appearing here is
$$
\sum_{j=0}^{r-1} (-1)^j {r-1 \choose j} z^{r-1-j} t_r^j \sum_{k_rt_r\le z-t_1-\cdots -t_{r-1} \atop{(k_r,a_r)=1 \atop{k_r\ge 1}}} k_r^j
$$
$$
= \sum_{j=0}^{r-1} (-1)^j {r-1 \choose j} z^{r-1-j} t_r^j \left( \frac{(z-t_1-\cdots -t_{r-1})^{j+1}\phi(a_r)}{(j+1)t_r^{j+1} a_r} +O(z^j \theta(a_r))\right)
$$
$$
=\frac{\phi(a_r)}{t_ra_r} z^r \sum_{j=0}^{r-1} (-1)^j {r-1 \choose j} \frac1{j+1} + O(  z^{r-1}\theta(a_r))
$$
$$
=\frac{\phi(a_r)}{rt_ra_r} z^r + O(z^{r-1}\theta(a_r))
$$
and the proof is complete.

\vskip2mm
{\bf Proof of Theorem 1.}
Apply Lemma 1 for $z=\log x$, $t_1=\log p_1,...,t_r=\log p_r$.
\vskip2mm

In order to prove Theorem 2 we need

\vskip2mm
{\bf Lemma 2.} {\it Let $r\ge 1$ be a fixed integer and $t_1,...,t_r>0$ be fixed
real numbers. Then for $z\ge 3$,
$$
\sum_{k_1t_1+\cdots +k_rt_r\le z \atop{k_1,...,k_r\ge 1}}
\prod_{i=1}^r \frac{\phi(k_i)}{k_i} =
\frac1{r! (\zeta(2))^r} \left(\prod_{i=1}^r \frac1{t_i}\right) z^r + O\left( z^{r-1}\log z\right). \leqno(6)
$$ }
\vskip2mm

\vskip2mm {\bf Proof of Lemma 2.} Induction on $r$, similar to the proof of Lemma 1.
We use the well-known estimate: let $s\ge -1$ be a real number, then for $z\ge 3$,
$$
\sum_{n\le z} \phi(n)n^s =\frac{z^{s+2}}{(s+2)\zeta(2)} + O(z^{s+1}\log z). \leqno(7)
$$

For $r=1$ (6) follows from (7) applied for $s=-1$. Suppose formula (6)
is valid for $r-1$ and prove it for $r$.
$$
\sum_{k_1t_1+\cdots +k_rt_r\le z \atop{k_1,...,k_r\ge 1}} \prod_{i=1}^r \frac{\phi(k_i)}{k_i} = \sum_{k_rt_r\le z-t_1-\cdots -t_{r-1} \atop{k_r\ge 1}} \frac{\phi(k_r)}{k_r}
\sum_{k_1t_1+\cdots +k_{r-1}t_{r-1}\le z-k_rt_r \atop{k_1,...,k_{r-1}\ge 1}} \prod_{i=1}^{r-1} \frac{\phi(k_i)}{k_i}
$$
$$
=\sum_{k_rt_r\le z-t_1-\cdots -t_{r-1} \atop{k_r\ge 1}} \frac{\phi(k_r)}{k_r}
\left(\frac1{(r-1)! (\zeta(2))^{r-1}} \left( \prod_{i=1}^{r-1} \frac1{t_i} \right)
(z-k_rt_r)^{r-1} + O\left( z^{r-2} \log z\right) \right)
$$
$$
= \frac1{(r-1)! (\zeta(2))^{r-1}} \left( \prod_{i=1}^{r-1} \frac1{t_i}\right)
\sum_{k_rt_r\le z-t_1-\cdots -t_{r-1} \atop{k_r\ge 1}} \frac{\phi(k_r)}{k_r}
(z-k_rt_r)^{r-1}+ O\left( z^{r-1} \log z \right).
$$
The sum appearing here is, applying (7),
$$
\sum_{j=0}^{r-1} (-1)^j {r-1 \choose j} z^{r-1-j} t_r^j \sum_{k_rt_r\le z-t_1-\cdots
-t_{r-1} \atop{k_r\ge 1}} \phi(k_r)k_r^{j-1}
$$
$$
= \sum_{j=0}^{r-1} (-1)^j {r-1 \choose j} z^{r-1-j} t_r^j \left( \frac{(z-t_1-\cdots
 -t_{r-1})^{j+1}}{(j+1)t_r^{j+1} \zeta(2)} +O(z^j \log z)\right)
$$
$$
=\frac1{t_r\zeta(2)} z^r \sum_{j=0}^{r-1} (-1)^j {r-1 \choose j} \frac1{j+1} + O(z^{r-1}\log z)
$$
$$
=\frac1{r t_r \zeta(2)} z^r + O(z^{r-1}\log z),
$$
which completes the proof.
\vskip2mm

\vskip2mm
{\bf Lemma 3.} {\it Let $r\ge 1$ be a fixed integer and $t_1,...,t_r>0$ be fixed
real numbers. Then for $z\ge 3$,
$$
\sum_{k_1t_1+\cdots +k_rt_r\le z \atop{j_1t_1+\cdots +j_rt_r\le z \atop{(k_1,j_1)=\cdots
=(k_r,j_r)=1 \atop{k_1,j_1,...,k_r,j_r\ge 1}}}} 1=
\frac1{(r!)^2 (\zeta(2))^r} \left(\prod_{i=1}^r \frac1{t_i^2}\right) z^{2r} +
O\left( z^{2r-1}\log z\right). \leqno(8)
$$ }
\vskip2mm

\vskip2mm
{\bf Proof of Lemma 3.} Using estimate (4),
$$
\sum_{k_1t_1+\cdots +k_rt_r\le z \atop{j_1t_1+\cdots +j_rt_r\le z \atop{(k_1,j_1)=\cdots =(k_r,j_r)=1 \atop{k_1,j_1,...,k_r,j_r\ge 1}}}} 1
=
\sum_{k_1t_1+\cdots +k_rt_r\le z \atop{k_1,...,k_r\ge 1}}
\sum_{j_1t_1+\cdots +j_rt_r\le z \atop{(j_1,k_1)=\cdots =(j_r,k_r)=1 \atop{j_1,...,j_r\ge 1}}} 1
$$
$$
=\sum_{k_1t_1+\cdots +k_rt_r\le z \atop{k_1,...,k_r\ge 1}}
\left(\frac1{r!} \left( \prod_{i=1}^r \frac{\phi(k_i)}{k_i t_i} \right) z^r +
O\left(z^{r-1} \sum_{i=1}^r \theta(k_i) \right)\right)
$$
$$
=\frac{z^r}{r!\prod_{i=1}^r t_i} \sum_{k_1t_1+\cdots +k_rt_r\le z \atop{k_1,...,k_r\ge 1}}
\prod_{i=1}^r \frac{\phi(k_i)}{k_i} + O\left( z^{r-1} \sum_{i=1}^r \sum_{k_1t_1+\cdots +k_rt_r\le z} \theta(k_i) \right)
$$
here the $O$-term is $O(z^{r-1}z^{r-1}z\log z)=O(z^{2r-1}\log z)$ and applying Lemma
2 to the main term finishes the proof. \vskip2mm

\vskip2mm
{\bf Proof of Theorem 2.} Apply Lemma 3 for $z=\log x$, $t_1=\log p_1,...,t_r=\log p_r$.
\vskip2mm

\vskip2mm {\bf Proof of Corollary 2.} This is a direct consequence of Theorem 2 and
Corollary 1. The considered asymptotic density is
$$
\lim_{x\to \infty} P^{(e)}(p_1,...,p_r;x) (N(p_1,...,p_r;x))^{-2}=
(\zeta(2))^{-r}.
$$

\vskip4mm

\noindent{{\bf L\'aszl\'o T\'oth}\\
University of P\'ecs\\
Institute of Mathematics and Informatics\\
Ifj\'us\'ag u. 6\\
7624 P\'ecs, Hungary\\
ltoth@ttk.pte.hu}

\end{document}